\documentclass[11pt]{epiarticle}
\usepackage{epimath}
\usepackage{cancel}
\usepackage{pgf,tikz}
\usetikzlibrary{shapes,arrows,shadows,patterns,shapes.geometric,calc}
\usepackage{graphicx}

\setpapertype{A4}

\usepackage[english]{babel}

\newtheorem{hypothesis}{Hypothesis}

\title{\LARGE{A \textit{pictorial} proof of the Four Colour Theorem}}

\titlemark{Four Colour Theorem}
\author{\normalsize{Bhupinder Singh Anand\footnote{\# 1003, Lady Ratan Tower, Dainik Shivner Marg, Gandhinagar, Worli, Mumbai - 400 018, Maharashtra, India. Email: bhup.anand@gmail.com. Mbl: +91 93225 91328. Tel: +91 (22) 2491 9821.}}}
\authormark{B S Anand}
\date{\tiny{\href{https://orcid.org/0000-0003-4290-9549}{https://orcid.org/0000-0003-4290-9549}}} 
\journal{Update as on \today.} 

\begin{document}
	
	\maketitle
	
	
	\begin{prelims}
		
		\def\abstractname{Abstract}
		\abstract{We give a \textit{pictorial}, and absurdly simple, proof that transparently illustrates \textit{why} four colours suffice to chromatically differentiate any set of contiguous, simply connected and bounded, planar spaces; by showing that there is no \textit{minimal} planar map. We show, moreover, why the proof cannot be expressed within classical graph theory.} \\  \\ \MSCclass 05C15  \\ \tiny \textbf{DECLARATIONS} $\bullet$ \textbf{Funding}: Not applicable $\bullet$ \textbf{Conflicts of interest/Competing interests}: Not applicable $\bullet$ \textbf{Availability of data and material}: Not applicable $\bullet$ \textbf{Code availability}: Not applicable $\bullet$ \textbf{Authors' contributions}: Not applicable
		
		
	\end{prelims}
	

	\section{Introduction}
	\label{sec:intro.fct}
	
	Although the Four Colour Theorem 4CT is considered pass\'{e} (see \S \ref{sec:fct.hist.prspctv}), we give a \textit{pictorial}, and absurdly simple, proof\footnote{Extracted from \cite{An21}, \S 1.G: \textit{\textit{Evidence-based} (pictorial), \textit{pre-formal}, proofs of the Four Colour Theorem}.} that transparently illustrates \textit{why} four colours suffice to chromatically differentiate any set of contiguous, simply connected and bounded, planar spaces; by showing that:
	
	\begin{enumerate}
		\item[(1)] If, for some natural numbers $m, n$, every planar map of less than $m + n$ contiguous, simply connected and bounded, areas can be $4$-coloured;
		
		\item[(2)] And, we assume (Hypothesis \ref{hyp:minimality}) that there is a \textit{sub-minimal} $4$-coloured planar map $\mathcal{M}$, of $m + n$ such areas, where \textit{finitary} creation of a \textit{specific}, additional, contiguous, simply connected and bounded, area $C$ within $\mathcal{M}$ yields a \textit{minimal} map $\mathcal{H}$ which entails that $C$ require a 5$^{th}$ colour;
		
		\item[(3)] Then Hypothesis \ref{hyp:minimality} is false (by Theorem \ref{thm:4CT}), since there can be no such \textit{sub-minimal} $4$-coloured planar map $\mathcal{M}$.
	\end{enumerate}

	Moreover we show why---challenging deep-seated dogmas that seemingly yet await, even if not actively seek, a mathematically `insightful', and philosophically `satisfying', proof of 4CT \textit{within} inherited paradigms---the \textit{pictorial} proof cannot be expressed within classical graph theory.

	\subsection{A historical perspective}
	\label{sec:fct.hist.prspctv}
	
	It would probably be a fair assessment that the mathematical significance of any new proof of the Four Colour Theorem 4CT continues to be perceived as lying not in any ensuing theoretical or practical utility of the Theorem per se, but in whether the proof can address the philosophically `unsatisfying', and occasionally `despairing' (see \cite{Tym79}; \cite{Sw80}; \cite{Gnt08}, \cite{Cl01}), lack of mathematical `insight', `simplicity', and `elegance' in currently known proofs of the Theorem (eg. \cite{AH77}, \cite{AHK77}, \cite{RSST}, \cite{Gnt08})---an insight and simplicity this investigation seeks in a \textit{pre-formal}\footnote{The need for distinguishing between \textit{belief-based} `informal', and \textit{evidence-based} `pre-formal', reasoning is addressed by Markus Pantsar in \cite{Pan09}; see also \cite{An21}, \S 1.D.} proof of 4CT.   
	
	\vspace{+1ex}
	For instance we note---amongst others---some candid comments from Robertson, Sanders, Seymour, and Thomas's 1995-dated (apparently pre-publication) summary\footnote{See \cite{RSSp}; also \cite{Thm98}, \cite{Cl01}, and the survey \cite{Rgrs} by Leo Rogers.} of their proof \cite{RSST}:
	
	\begin{quote}
		\footnotesize
		``\textbf{Why a new proof?}
		
		\vspace{+1ex}
		There are two reasons why the Appel-Haken proof is not completely satisfactory.
		
		\begin{quote}
			$\bullet$ Part of the Appel-Haken proof uses a computer, and cannot be verified by hand, and
			
			$\bullet$ even the part that is supposedly hand-checkable is extraordinarily complicated and tedious, and as far as we know, no one has verified it in its entirety." \textit{\tiny{\ldots Robertson et al: \cite{RSSp}, Pre-publication.}}
		\end{quote}	
		
		``It has been known since 1913 that every minimal counterexample to the Four Color Theorem is an internally six-connected triangulation. In the second part of the proof, published in [\textbf{4}, p. 432], Robertson et al. proved that at least one of the 633 configurations appears in every internally six-connected planar triangulation. This condition is called ``unavoidability," and uses the discharging method, first suggested by Heesch. Here, the proof differs from that of Appel and Haken in that it relies far less on computer calculation. Nevertheless, parts of the proof still cannot be verified by a human. The search continues for a computer-free proof of the Four Color Theorem." \textit{\tiny{\ldots Brun: \cite{Bru02}, \S 1. Introduction (Article for undergraduates)}}
		
		\vspace{+1ex}
		``The four-colour problem had a long life before it eventually became the four-colour theorem. In 1852 Francis Gutherie (later Professor of Mathematics at the University of Cape Town) noticed that a map of the counties of England could be coloured using only four colours. He wondered if four colours would always suffice for any map. He, or his brother Frederick, proposed the problem to Augustus De Morgan (see the box at the end of Section 3.5 in Chapter 3) who liked it and suggested it to other mathematicians. Interest in the problem increased after Arthur Cayley presented it to the London Mathematical Society in 1878 (\cite{Cay79}). The next year Alfred Bray Kempe (a British lawyer) gave a proof of the conjecture. His proof models the problem in terms of graphs and breaks it up into a number of necessary cases to be checked. Another proof was given by Peter Tait in 1880. It seemed that the four-colour problem had been settled in the affirmative.
		
		\vspace{+1ex}
		However, in 1890 Percy John Heawood found that Kempe's proof missed one crucial case, but that the approach could still be used to prove that five colours are sufficient to colour any map. In the following year Tait's proof was also shown to be flawed, this time by Julius Petersen, after whom the Petersen graph is named. The four-colour problem was therefore again open, and would remain so for the next 86 years. In that time it attracted a lot of attention from professional mathematicians and good (and not so good) amateurs alike. In the words of Underwood Dudley:
		
		\begin{quote}
			The four-color conjecture was easy to state and easy to understand, no large amount of technical mathematics is needed to attack it, and errors in proposed proofs are hard to see, even for professionals; what an ideal combination to attract cranks!
		\end{quote}
		
		The four-colour theorem was finally proved in 1976 by Kenneth Appel and Wolfgang Haken at the University of Illinois. They reduced the problem to a large number of cases, which were then checked by computer. This was the first mathematical proof that needed computer assistance. In 1997 N. Robertson, D.P. Sanders, P.D. Seymour and R. Thomas published a refinement of Appel [and] Haken's proof, which reduces the number of necessary cases, but which still relies on computer assistance. The search is still on for a short proof that does not require a computer." \textit{\tiny{\ldots Conradie/Goranko: \cite{CG15}, \S 7.7.1, Graph Colourings, p.417.}}
		
		\vspace{+1ex}
		``The first example concerns a notorious problem within the philosophy of mathematics, namely the acceptability of computer-generated proofs or proofs that can only be checked by a computer; for instance because it includes the verification of an excessively large set of cases. The text-book example of such a mathematical result is the proof of the 4-colour theorem, which continues to preoccupy philosophers of mathematics (Calude 2001). Here, we only need to note that the debate does not primarily concern the correctness of the result, but rather its failure to adhere to the standard of \textit{surveyability} to which mathematical proofs should conform." \textit{\tiny{\ldots Allo: \cite{All17}, Conclusion, p.562.}}
		
		\vspace{+1ex}
		``Being the first ever proof to be achieved with substantial help of a computer, it has raised questions to what a proof really is. Many mathematicians remain sceptical about the nature of this proof due to the involvement of a computer. With the possibility of a computing error, they do not feel comfortable relying on a machine to do their work as they would be if it were a simple pen-and-paper proof.
		
		\vspace{+1ex}
		The controversy lies not so much on whether or not the proof is valid but rather whether the proof is a valid proof. To mathematicians, it is as important to understand why something is correct as it is finding the solution. They hate that there is no way of knowing how a computer reasons. Since a computer runs programs as they are fed into it, designed to tackle a problem in a particular way, it is likely they will return what the programmer wants to find leaving out any other possible outcomes outside the bracket.
		
		\vspace{+1ex}
		Many mathematicians continue to search for a better proof to the problem. They prefer to think that the Four Colour problem has not been solved and that one day someone will come up with a simple completely hand checkable proof to the problem." \textit{\tiny{\ldots Nanjwenge: \cite{Nnj18}, Chapter 8, Discussion (Student Thesis).}}
		
		\vspace{+1ex}
		``The heavy reliance on computers in Appel and Haken's proof was immediately a topic of discussion and concern in the mathematical community. The issue was the fact that no individual could check the proof; of special concern was the reductibility [\textit{sic}] part of the proof because the details were ``hidden" inside the computer. Though it isn't so much the \textit{validity} of the result, but the \textit{understanding} of the proof. Appel himself commented: ``\ldots there were people who said, `This is terrible mathematics, because mathematics should be clean and elegant,' and I would agree. It would be nicer to have clean and elegant proofs." See page 222 of Wilson." \textit{\tiny{\ldots Gardner: \cite{Grd21}, \S 11.1, Colourings of Planar Maps, pp.6-7 (Lecture notes).}}
	\end{quote}

	\section{A \textit{pictorial} proof of the 4-Colour Theorem}
	
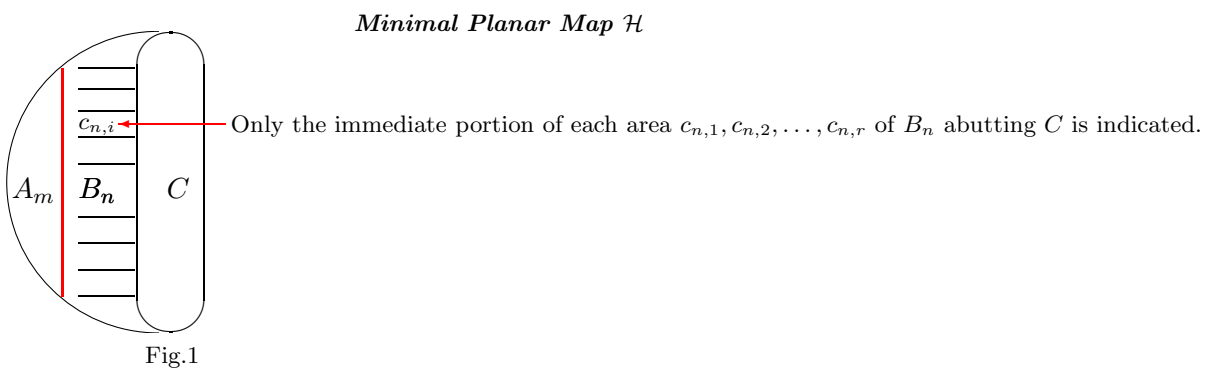
\label{sec:prf.fct}

	\begin{picture}(50,20)
		
		\begin{tikzpicture}[baseline=(current bounding box.north)]
			
			\begin{scope}
				\clip (-2,-2) rectangle (0,2);
				\draw (0,0) circle(2);
				\draw (-2,10) -- (2,10);
			\end{scope}
			%
		\end{tikzpicture}
		
		\put(4.45,-57){\oval(25,113)}
		\put(-36,-100){\color{red}\line(0,1){86}}
		
		\put(-30,-100){\line(1,0){21}}
		\put(-30,-90){\line(1,0){21}}
		\put(-30,-80){\line(1,0){21}}
		\put(-30,-70){\line(1,0){21}}
		\put(-30,-50){\line(1,0){21}}
		\put(-30,-40){\line(1,0){21}}
		\put(-30,-30){\line(1,0){21}}
		\put(-30,-22){\line(1,0){21}}
		\put(-30,-14){\line(1,0){21}}
		
		\put(25,-35){\color{red}\vector(-1,0){40}}
		
		\put(3,-63){$C$}
		\put(-30,-63){$B_{n}$}
		\put(-55,-63){$A_{m}$}
			
		\put(73,0){\textbf{\textit{\footnotesize Minimal Planar Map $\mathcal{H}$}}}
		
		\put(-30,-63){$B_{n}$}
		
		\put(27,-38){\footnotesize{Only the immediate portion of each area $c_{n, 1}, c_{n,2}, \ldots, c_{n, r}$ of $B_n$ abutting $C$ is indicated.}}
		\put(-30,-37){\footnotesize $c_{n, i}$}
		
		\put(-5,-125){\footnotesize{Fig.1}}
	\end{picture}
	
\vspace{+30ex}	
\noindent We consider the surface of the hemisphere (\textit{minimal} planar map $\mathcal{H}$) in Fig.1 where:
		
\begin{enumerate}
	\item $A_m$ denotes a region of $m$ contiguous, simply connected and bounded, surface areas $a_{m,1}, a_{m, 2},$ $\ldots, a_{m, m}$, \textit{none} of which share a non-zero boundary segment with the contiguous, simply connected, surface area $C$ (as indicated by the red barrier which, however, is \textit{not} to be treated as a boundary of the region $A_m$);
			
	\item $B_n$ denotes a region of $n$ contiguous, simply connected and bounded, surface areas $b_{n,1}, b_{n, 2},$ $\ldots, b_{n, n}$, \textit{some} of which, say $c_{n, 1}, c_{n, 2}, \ldots, c_{n, r}$, share \textit{at least one} non-zero boundary segment of $c_{n, i}$ with $C$; where, for each $1 \leq i \leq r$, $c_{n, i} = b_{n, j}$ for some $1 \leq j \leq n$;
			
	\item $C$ is a single contiguous, simply connected and bounded, area created \textit{finitarily} (i.e., not \textit{postulated}) by sub-dividing and annexing one or more contiguous, simply connected, portions of each area $c_{n, i}^-$ (defined in Hypothesis \ref{hyp:minimality}(b) below) in the region $B_n^-$ (defined in Hypothesis \ref{hyp:minimality}(b)).
\end{enumerate}

\begin{hypothesis}[Minimality Hypothesis]
	\label{hyp:minimality}
	Since four colours suffice for maps with fewer than 5 regions, we assume the existence of some $m, n$, in a putatively \textit{minimal} planar map $\mathcal{H}$, which defines a \textit{minimal} configuration of the region $\{A_m + B_n + C\}$ where:
	
	\begin{enumerate}		
		\item[(a)] \textit{any} configuration of $p$ contiguous, simply connected and bounded, areas \textit{can} be $4$-coloured if $p \leq m + n$, where $p, m, n \in \mathbb{N}$, and $m+n \geq 5$;
		
		\item[(b)] \textit{any} configuration of the $m + n$ contiguous, simply connected and bounded, areas of the region, say $\{A_m^- + B_n^-\}$, in a putative, \textit{sub-minimal}, planar map $\mathcal{M}$ \textit{\underline{before}} the creation of $C$---constructed \textit{finitarily} by sub-dividing and annexing some portions from each area, say $c_{n, i}^-$, of $B_n^-$ in $\mathcal{M}$---\textit{can} be 4-coloured;
		
		\item[(c)] the region $\{A_m + B_n +C\}$ in the planar map $\mathcal{H}$ is a \textit{\underline{specific}} configuration of $m + n + 1$ contiguous, simply connected and bounded, areas that \textit{cannot} be 4-coloured (whence the area $C$ \textit{necessarily} requires a 5$^{th}$ colour by the \textit{Minimality Hypothesis}).	
	\end{enumerate}
\end{hypothesis}

\begin{theorem}[Four Colour Theorem]	
	\label{thm:4CT}
	\textit{No planar map needs more than four colours.}
\end{theorem}

\begin{picture}(50,20)
	
	\begin{tikzpicture}[baseline=(current bounding box.north)]
		
		\begin{scope}
			\clip (-2,-2) rectangle (0,2);
			\draw (0,0) circle(2);
			\draw (-2,10) -- (2,10);
		\end{scope}
		%
	\end{tikzpicture}
	
	\put(4.45,-57){\oval(25,113)}
	\put(-36,-100){\color{red}\line(0,1){86}}
	
	\put(-30,-100){\line(1,0){21}}
	\put(-30,-90){\line(1,0){21}}
	\put(-30,-80){\line(1,0){21}}
	\put(-30,-70){\line(1,0){21}}
	\put(-30,-50){\line(1,0){21}}
	\put(-30,-40){\line(1,0){21}}
	\put(-30,-30){\line(1,0){21}}
	\put(-30,-22){\line(1,0){21}}
	\put(-30,-14){\line(1,0){21}}

	\put(3,-17){\tiny $E_1$}
	\put(3,-93){\tiny $E_1$}
	\put(-5,-37){\tiny $D_1$}

	\put(3,-63){\tiny $E_1$}
	\put(-30,-63){$B_{n}$}
	\put(-55,-63){$A_{m}$}

	\put(73,0){\textbf{\textit{\footnotesize Planar Map $\mathcal{H'}$}}}
	
	\put(-30,-63){$B_{n}$}

	\put(-30,-37){\footnotesize $c_{n, i}$}
	\put(-10,-40){\line(1,0){16}}
	\put(-10,-30){\line(1,0){16}}
	\put(5,-40){\line(0,1){10}}

	\put(-5,-125){\footnotesize{Fig.2}}
\end{picture}

\vspace{+27ex}

\begin{proof}
	If the area $C$ of the \textit{minimal} planar map $\mathcal{H}$ in Fig.1 is divided further (as indicated in Fig.2) into two \textit{non-empty} areas $\small D_1$ and $\small E_1$, where:
	
	\begin{itemize}
		\item $\small D_1$ shares a non-zero boundary segment with only \textit{one} of the areas $c_{n, i}$; and
		
		\item $D_1$ can be \textit{treated} as an original area of $c_{n, i}^-$ in $\mathcal{M}$ (see Hypothesis \ref{hyp:minimality}(b)) that was annexed to form part of $C$ in $\mathcal{H}$ (in Fig.1);
	 \end{itemize}
	 
	\noindent then $D_1$ can be absorbed back into $c_{n, i}$ without violating the Minimality Hypothesis. Moreover, $c_{n, i} + \small{D_1}$ \textit{must} share a non-zero boundary with $\small E_1$ in $\mathcal{H'}$ if $c_{n, i} = b_{n,j}$ for some $1<j<n$, and $b_{n,j}, C$ are required to be differently coloured, in $\mathcal{H}$.
	
	\vspace{+1ex}
	Such a division, as illustrated in Fig.2, \textit{followed} by re-absorption of $\small D_1$ into $c_{n, i}^-$ (denoted, say, by $B_{n} + \small{D_1}$), would reduce the configuration $\mathcal{H'}$ in Fig.2 again to a \textit{minimal} planar map, say $\mathcal{H}_1$ with a configuration $\{A_m + B_{n}' + \small{E_1}\}$, where $B_{n}' = (B_{n} + \small{D_1})$; which would in turn necessitate a 5$^{th}$ colour for the area $\small E_1 \subset C$ by the Minimality Hypothesis.
	
	\vspace{+1ex}
	Since we cannot, by reiteration, have a non-terminating sequence $C \supset \small{E_1} \supset \small{E_2} \supset \small{E_3} \supset \ldots$, the sequence must terminate in a \textit{non-empty} area $\small{E_k}$ of a \textit{minimal} planar map, say $\mathcal{H}_k$, for some finite integer $k$; where $\small E_k$ contains no area that is annexed from any of the areas of $B_n^-$ in $\mathcal{M}$ prior to the formation of the \textit{minimal} planar map $\mathcal{H}$ (in Fig.1).
	
	\begin{quote}
		\footnotesize
		\textbf{Comment}: Note that we cannot admit as a putative limit of $C \supset \small{E_1} \supset \small{E_2} \supset \small{E_3} \ldots$ the configuration where all the $c^{-}_{n, i}$---corresponding to the abutting areas $c_{n, i}$ of $C$ in the Minimal Planar Map $\mathcal{H}$---meet at a point in the putative, sub-minimal, planar map $\mathcal{M}$, since any \textit{finitary} (i.e., not postulated) creation of $C$, begun by initially annexing a non-empty area of some $c^{-}_{n, i}$ at such an apex (corresponding to the putative `finally merged' area of the above non-terminating sequence $C \supset \small{E_1} \supset \small{E_2} \supset \small{E_3} \supset \ldots$), would require, at most, a 4$^{th}$ but not a 5$^{th}$ colour.
	\end{quote}
	
	However, by Hypothesis \ref{hyp:minimality}(b), this contradicts the definition of the area $C$, in the \textit{minimal} planar map $\mathcal{H}$ (in Fig.1)---ergo of the area $\small{E_k}$ in the \textit{minimal} planar map $\mathcal{H}_k$---as formed \textit{finitarily} by sub-division and annexation of existing areas of $B_n^-$ in $\mathcal{M}$.
	
	\vspace{+1ex}
	Hence there can be no \textit{minimal} planar map $\mathcal{H}$ which defines a \textit{minimal} configuration such as the region $\{A_m + B_n + C\}$ in Fig.1. The theorem follows. \hfill $\Box$
\end{proof}

We conclude by noting that, since classical graph theory\footnote{See, for instance, Brun \cite{Bru02}, Conradie/Goranko \cite{CG15}, Gardner \cite{Grd21}.} represents non-empty areas as points (\textit{vertices}), and a non-zero boundary between two areas as a line (\textit{edge}) joining two points (\textit{vertices}), it cannot express the proof of Theorem \ref{thm:4CT} graphically.

\vspace{+1ex}
\textit{Reason}: The proof of Theorem \ref{thm:4CT} appeals to finitarily distinguishable properties of a series of, putatively \textit{minimal}, planar maps $\mathcal{H}, \mathcal{H}_1, \mathcal{H}_2, \ldots$, created by a corresponding sequence of areas $C, \small{E_1}, \small{E_2}, \ldots,$ where each area is finitarily created from, or as a proper subset of, some preceding area/s in such a way that the graphs of $\mathcal{H}, \mathcal{H}_1, \mathcal{H}_2, \ldots$ remain undistinguished.

\end{document}